\magnification=\magstephalf
\baselineskip=18truept
\abovedisplayskip=6pt plus 3pt minus 3pt
\belowdisplayskip=6pt plus 3pt minus 3pt
\overfullrule=0pt
\def\n{{\bf n}}
\def\x{{\bf x}}
\def\m{{\bf m}}
\def\g{\gamma}
\def\o{\over}
\def\a{\alpha}
\def\lm{\lambda}
\def\ss{\scriptstyle}
\centerline{\bf Some Systems of Multivariable Orthogonal}
\centerline{\bf q-Racah polynomials}
\medskip
\centerline{George Gasper\footnote{*}{Department of Mathematics, 
Northwestern University,
Evanston, IL 60208.} and Mizan Rahman \footnote{\dag}{School of 
Mathematics and Statistics,
Carleton University, Ottawa, ON, K1S 5B6, CANADA. Supported in part
by NSERC grant \#A6197.}}
\medskip
\noindent
{\bf Abstract}.\quad In 1991 Tratnik derived two systems of 
multivariable orthogonal
Racah polynomials and considered their limit cases. $q$-Extensions of 
these systems
are derived, yielding systems of multivariable orthogonal $q$-Racah
polynomials,  from which systems of multivariable orthogonal\break 
$q$-Hahn, dual $q$-Hahn,
$q$-Krawtchouk,
$q$-Meixner, and $q$-Charlier polynomials follow as special or limit
cases.
$q$-Hahn, 
$q$-Krawtchouk, $q$-Meixner,
and $q$-Charlier polynomials
\medskip
\bigskip
\noindent
{\bf 1. Introduction}

\noindent
The Racah polynomials [1, 2, 25], defined by
$$\eqalignno{
r_n(x) &= r_n(x; \a, \beta, \g, N)\cr
&= (\a+1)_n (\beta + \g + 1)_n (-N)_n\cr
&\quad \times \, _4F_3 \left[\matrix{-n, n + \a + \beta +1, -x, x + \g -N\cr
\a+1, \beta + \g + 1, -N\cr} ;1\right]&(1.1)\cr
}
$$
for $n = 0, 1, \ldots, N$, satisfy the discrete orthogonality relation
$$\sum^N_{x=0} r_n (x) r_m (x) \rho (x) = \lm_n \delta_{n,m}\eqno (1.2)$$
for $n, m = 0, 1, \ldots, N$, with the weight function
$$\eqalignno{
\rho (x) &= \rho (x; \a, \beta, \g, N)\cr
&= {\g -N+2x\o \g -N}\ \, {(\g - N)_x (\a + 1)_x (\beta + \g +1)_x (-N)_x\o
x!\,(\g - \a - N)_x (-\beta - N)_x (\g + 1)_x} &(1.3)\cr
}
$$
and the normalization constant
$$\eqalignno{
\lm_n &= \lm _n (\a, \beta, \g, N)\cr
&= {(\a + \beta + 2)_N (-\g)_N\o
(\beta + 1)_N (\a -\g + 1)_N} (\a -\g + 1)_n(\beta + \g + 1)_n\, (\a 
+ \beta + 2 + N)_n\cr
&\quad \times {(\a + \beta + 1)n!\,(\a + 1)_n (\beta +1)_n (-N)_n\o
(\a + \beta + 1 + 2n) (\a + \beta + 1)_n}, &(1.4)\cr
}
$$
where $N$ is a nonnegative integer. In 1991 Tratnik [24] extended the 
Racah polynomials to a system of multivariable orthogonal polynomials 
(in a slightly different notation)
$$\eqalignno{
&R_{\n} (\x) = R_\n (\x; a_1, \ldots, a_{s+1}, \eta, N)\cr
&= \prod^s_{k=1} r_{n_k} (x_k - N_{k-1}; 2N_{k-1}+ \eta + \a_k - a_1, 
\a_{k+1} -1, N_{k-1} + \a_k + x_{k+1}, x_{k+1} - N_{k-1}),&(1.5)\cr
}
$$
where
$$\eqalignno{
&\x = (x_1, \ldots, x_s), \ x_{s+1} = N,\ \n = (n_1, \ldots, n_s), N_0 = 0,\cr
&N_k = \sum^k_{j=1} n_j, \ 1 \le k \le s,\ \ \ \ \a_k = \sum^k_{j=1} 
a_j, \  1 \le k \le s+1,&(1.6)\cr
}
$$
and $N_s \le N$. Clearly, $R_{\n}(\x)$ is a polynomial of total 
degree $N_s$ in the variables $y_1, \ldots, y_s$ with $y_k = x_k (x_k 
+ \a_k)$, $k = 1, \ldots, s$.
\medskip
Tratnik showed that these polynomials satisfy the discrete 
orthogonality relation
$$\sum_{\x} R_{\n} (\x) R_{\m} (\x) \rho (\x) = \lm_n 
\delta_{\n,\m}\eqno (1.7)$$
for $N_s, M_s \le N$, where $M_k = \sum^k_{j=1} m_j$ and the summation is over
  all $\x = (x_1, \ldots, x_s)$ with $x_k = 0, 1, \ldots, N$ for $k = 
1, \ldots, s$, $\delta_{\n, \m} =
\prod\limits^s_{k=1} \delta_{n_k, m_k}$, the weight function is
$$\eqalignno{
\rho (\x) &= \rho (\x; a_1, \ldots, a_{s+1}, \eta, N)\cr
&= {N!\,\Gamma (\a_s + N+1)\o
\Gamma (a_{s+1} + N) \Gamma (\a_{s+1} + N)}\, {(a_1)_{x_1} (\eta + 1)_{x_1}\o
x_1 !\,(a_1 -\eta)_{x_1}}\cr
&\quad \times \prod^s_{k=1} {\Gamma (a_{k+1} + x_{k+1} - x_k) \Gamma (\a_{k+1} 
+ x_{k+1} + x_k)\o
(x_{k+1} - x_k)!\,\Gamma (\a_k + x_{k+1}+ x_k + 1)} {\a_k + 2x_k\o
\a_k}&(1.8)\cr
}
$$
and
$$\eqalignno{
\lm_n & =\lm_n (a_1, \ldots, a_{s+1}, \eta, N)\cr
&= (\a_s + N)_{N_s} (\eta + 1 -a_1 - N)_{N_s} (-N)_{N_s}\cr
&\quad \times {\Gamma (a_1 - \eta) \Gamma (\a_s + N+1) \Gamma 
(\a_{s+1} -a_1 + \eta + N_s + N+1)\o
\Gamma (a_1) \Gamma (\eta + 1) \Gamma (a_{s+1} + N) \Gamma (a_1 - \eta + N)}\cr
&\quad \times \left[ \prod^s_{k=1} \a^{-1}_k n_k!\,(\a_{k+1} - a_1 + 
\eta + N_k +
N_{k-1})_{n_k}\right.\cr &\qquad \times\left. {\Gamma (a_{k+1} + n_k) 
\Gamma (\a_k - a_1 + N_k +
N_{k-1}  + 1)\o
\Gamma (\a_{k+1} - a_1 + \eta + 2N_k + 1)}\right].&(1.9)\cr
}
$$
Note that $\rho (\x) = 0$ if $x_{k+1} < x_k$ for some $k < s$.
He also pointed out the special case of the multivariable Hahn 
polynomials of Karlin and McGregor [14],
the limit cases of the multivariable Krawtchouk, Meixner and Charlier 
polynomials, and used
permutations of the parameters and variables to derive a second 
system of polynomials that are
orthogonal with respect to the weight function in (1.8).
\medskip
In this paper we derive $q$-extensions of Tratnik's $R_{\n} (\x; a_1, 
\ldots, a_{s+1}, \eta, N)$ polynomials and of his second system of 
multivariable orthogonal Racah polynomials, from which systems of 
multivariable orthogonal $q$-Hahn, dual $q$-Hahn, $q$-Krawtchouk, 
$q$-Meixner, and $q$-Charlier polynomials follow as special or limit 
cases. Some $q$-extensions of Tratnik's [23] multivariable orthogonal 
Wilson polynomials and of his [21] multivariable biorthogonal 
generalization of the Wilson polynomials are given in our papers [11] 
and [10], respectively.

The authors wish to thank Michael Schlosser for his useful comments on the original
version of this paper.
\bigskip
\noindent
{\bf 2. Multivariable orthogonal $q$-Racah polynomials}.
\medskip
\noindent
Analogous to (1.1) we define the $q$-Racah polynomials [2] by
$$\eqalignno{
r_n (x;q) &= r_n (x; a, b, c, N;q)\cr
&= (aq, bcq, q^{-N}; q)_n (q^N/c)^{n/2}\cr
&\quad \times\, _4\phi_3\left[\matrix{
q^{-n}, abq^{n+1}, q^{-x}, cq^{x-N}\cr
aq, bcq, q^{-N}\cr} ;q,q\right]&(2.1)\cr
}
$$
for $n = 0, 1, \ldots, N$, where
$$(a;q)_n = \prod^{n-1}_{k=0} (1-aq^k), \quad (a_1, \ldots, a_k; q)_n 
= \prod^k_{j=1} (a_j;q)_n,$$
and we use the notation of our book [9]. The power $(q^N/c)^{n/2}$ is 
chosen in (2.1) so that certain symmetry properties of the $q$-Racah 
polynomials are satisfied; for example, by the Sears transformation 
formula [9, (2.10.4)] it follows from (2.1) that
$$r_n (x;a, b, c, N;q) = r_n (N-x; b, a, c^{-1}, N;q).\eqno (2.2)$$
Note that $r_n (x;q)$ is a polynomials of degree $n$ in the variable 
$z = q^{-x} + cq^{x-N}$.
\medskip
Askey and Wilson [2] showed that the $q$-Racah polynomials satisfy 
the orthogonality relation [9, (7.2.18)]
$$\sum^N_{x=0} r_n (x;q) r_m (x;q) \rho (x;q) = \lm_n(q) 
\delta_{n,m}\eqno (2.3)$$
for $n, m = 0, 1, \ldots, N$, with
$$\eqalignno{
\rho (x;q) &= \rho (x; a, b, c, N;q)\cr
&= {1-cq^{2x-N}\o
1-c q^{-N}} {(cq^{-N}, aq, bcq, q^{-N};q)_x\o
(q, ca^{-1}q^{-N}, b^{-1}q^{-N},cq;q)_x}(abq)^{-x}&(2.4)\cr
}
$$
and
$$\eqalignno{
\lm_n(q) &= \lm_n(a, b, c, N;q)\cr
&= {(c^{-1}, abq^2;q)_N\o
(aq/c, bq;q)_N} (q, aq, bq, aq/c, bcq, q^{-N};q)_n\cr
&\quad \times {1-abq\o
1-abq^{2n+1}} {(abq^{N+2};q)_n\o
(abq;q)_n}.&(2.5)\cr
}
$$
\indent
Analogous to (1.5) we define the multivariable $q$-Racah polynomials by
$$\eqalignno{
&R_{\n} (\x;q) = R_{\n} (\x; a_1, \ldots, a_{s+1}, b, N;q)\cr
&= \prod^s_{k=1} r_{n_k} (x_k - N_{k-1}; bA_kq^{2N_{k-1}}/a_1,
a_{k+1} q^{-1}, A_k q^{x_{k+1} + N_{k-1}}, x_{k+1} - 
N_{k-1};q),&(2.6)\cr
}
$$
where, in addition to the definitions of $\x$, $\n$, $x_{s+1}$ and 
$N_k$ given in (1.6), we let
$$A_0 = 1,\quad A_k = \prod^k_{j=1} a_j, \quad k = 1, \ldots, s+1.\eqno (2.7)$$
Note that $R_{\n} (\x;q)$ is a polynomial of total degree $N_s$ in 
the variables $z_1, \ldots, z_s$ where
$$z_k = q^{-x_k} + A_k q^{x_k}, \quad k = 1, \ldots, s.$$
\indent
In order to derive an orthogonality relation for these polynomials we 
start by considering the following $q$-extensions of the weight 
function $\rho (\x)$ in (1.8):
$$\eqalignno{
\rho (\x;q) &= \rho (\x;a_1, \ldots, a_{s+1}, b, N;q)\cr
&= {(q, qA_s;q)_N\o
(a_{s+1}, A_{s+1};q)_N} {(a_1, bq;q)_{x_1}\o
(q, a_1/b;q)_{x_1}}\cr
&\quad\times \prod^s_{k=1} {(a_{k+1}; q)_{x_{k+1}- x_k} 
(A_{k+1};q)_{x_{k+1}+ x_k} (1-A_kq^{2x_k})\o
(q;q)_{x_{k+1}-x_k} (qA_k;q)_{x_{k+1} + x_k} (1-A_k)} 
(c_{k,s})^{-x_k},&(2.8)\cr
}
$$
where $c_{1,s}, \ldots, c_{s,s}$ are to be determined (see (2.15)) so 
that the orthogonality relation
$$\sum^N_{x_s = 0} \sum^{x_s}_{x_{s-1} = 0} \ldots \sum^{x_3}_{x_2 = 
0} \sum^{x_2}_{x_1 = 0} R_{\n} (\x;q) R_{\m} (\x;q)\rho (\x;q) = 
\lm_{\n} (q) \delta_{\n,\m}\eqno (2.9)$$
holds for $N_s$, $M_s \le N$ and certain normalization constants 
$\lm_{\n}(q)$, which will be given in (2.14).
\medskip
The summation over $x_1$ in (2.9) can be evaluated via (2.3),
which is derived in [9, \S7.2] by means of the terminating
$_6\phi_5$ summation [9, (2,4,2)], to obtain 
 for $s \ge 2$ that
$$\eqalignno{
&\sum^{x_2}_{x_1 = 0} {(a_1, bq;q)_{x_1} (a_2;q)_{x_2-x_1} 
(A_2;q)_{x_2 + x_1} (1-a_1q^{2x_1})\o
(q, a_1/b;q)_{x_1} (q)_{x_2 - x_1} (a_1 q;q)_{x_2 + x_1} (1-a_1)} 
(c_{1,s})^{-x_1}\cr
&\quad \times r_{n_1} (x_1; b, a_2 q^{-1}, a_1 q^{x_2}, x_2;q) 
r_{m_1} (x_1;b, a_2q^{-1}, a_1q^{x_2}, x_2;q)\cr
&= {1-a_2b\o
1-a_2 bq^{2n_1}} {(bq, a_2, q, A_2, bA_2q/a_1;q)_{n_1}\o
(a_2b;q)_{n_1}} \Big({b\o a_1}\Big)^{n_1} q^{n^2_1}\cr
&\quad\times {(A_2q^{n_1}, b A_2q^{n_1+1}/a_1;q)_{x_2}\o
(q;q)_{x_2 - n_1} (a_1/b;q)_{x_2 - n_1}} (bq^{2n_1 + 
1})^{-x_2}\delta_{n_1, m_1}&(2.10)\cr
}
$$
provided we choose $c_{1,s} = bq$ for $s \ge 2$. In the single 
variable case $s=1$ it is clear from (2.3) and (2.4) that $c_{1,1} = 
bq$, and hence we conclude that $c_{1,s} = bq$ for $s \ge 1$.
\medskip
By doing the summations over $x_1, x_2, \ldots, x_j$ in (2.9) for 
small $j$ one is led to conjecture that after summing over $x_1, 
\ldots, x_j$ in (2.9) and setting $c_{k,s} = a_k$ for $2 \le k \le 
s-1$ we have
$$\eqalignno{
&\left[ \prod^j_{k=1} {(1-bA_{k+1}/a_1)\o
(1-bA_{k+1}q^{2N_k}/a_1)} {(q, a_{k+1};q)_{n_k} (bA_kq/a_1; q)_{N_k + 
N_{k-1}}\o
(b A_{k+1}/a_1;q)_{N_k + N_{k-1}}} \delta_{n_k, m_k}\right]\cr
&\quad \times {(A_{j+1}, qbA_{j+1}/a_1; q)_{N_j + x_{j+1}}\o
(q, a_1/b;q)_{x_{j+1}-N_j}} \Big({b\o a_1}\Big)^{N_j} q^{N^2_j} 
\Big({bA_j q^{2N_j+1}\o a_1}\Big)^{-x_{j+1}} &(2.11)\cr
}
$$
as the sum for $j = 1, 2, \ldots, s-1$.
\medskip
To prove (2.11) by induction on $j$, assume that $j \le s-2$ for $s 
\ge 3$, multiply (2.11) by the remaining $x_{j+1}$-dependent part of 
the weight function and polynomials, and then sum over $x_{j+1}$ to 
get

$$\eqalignno{
&\left[\prod^j_{k=1} {(1-bA_{k+1}/a_1)\o
(1-bA_{k+1}q^{2N_k}/a_1)} {(q, a_{k+1};q)_{n_k} (qbA_k/a_1;q)_{N_k + N_{k-1}}\o
(bA_{k+1}/a_1;q)_{N_k + N_{k-1}}}\delta_{n_k, m_k}\right]\cr
&\quad \times \Big({b\o a_1}\Big)^{N_j} q^{N^2_j} 
\sum^{x_{j+2}}_{x_{j+1}=0}{(A_{j+1}, qbA_{j+1}/a_1;q)_{N_j+ x_{j+1}}\o
(q, a_1/b;q)_{x_{j+1}-N_j}} \Big({bA_jq^{2N_j+1}\o a_1}\Big)^{-x_{j+1}}\cr
&\qquad\times {(a_{j+2};q)_{x_{j+2}-x_{j+1}} (A_{j+2};q)_{x_{j+2} + 
x_{j+1}} (1-A_{j+1}q^{2x_{j+1}})\o
(q;q)_{x_{j+2}  - x_{j+1}} (qA_{j+1};q)_{x_{j+2} + x_{j+1}} 
(1-A_{j+1})} (a_{j+1})^{-x_{j+1}}\cr
&\quad \times r_{n_{j+1}} (x_{j+1} - N_j; bA_{j+1}q^{2N_j}/a_1, 
a_{j+2}q^{-1}, A_{j+1}q^{N_j + x_{j+2}}, x_{j+2}-N_j;q)\cr
&\quad \times r_{m_{j+1}} (x_{j+1} - N_j; bA_{j+1}q^{2N_j}/a_1, 
a_{j+2}q^{-1}, A_{j+1}q^{N_j + x_{j+2}}, x_{j+2} - N_j;q).&(2.12)\cr
}
$$
Note that the summand above vanishes for $0 \le x_{j+1} < N_j$. 
Hence, the sum is effectively from $x_{j+1} = N_j$ to $x_{j+2}$. 
Replacing $x_{j+1} - N_j$ by $x$, say, the sum in (2.12) then becomes
$$\eqalignno{
&{(A_{j+1}, qbA_{j+1}/a_1;q)_{2N_j} (a_{j+2};q)_{x_{j+2} - 
N_j}(A_{j+2};q)_{x_{j+2} + N_j} (1-A_{j+1}q^{2N_j})\o
(q;q)_{x_{j+2} - N_j} (qA_{j+1};q)_{x_{j+2} + N_j} (1-A_{j+1})}\cr
&\times \Big({bA_{j+1}\o a_1} q^{2N_j + 1}\Big)^{-N_j} \sum^{x_{j+2} 
- N_j}_{x=0} {(1-A_{j+1} q^{2N_j + 2x})\o
(1-A_{j+1} q^{2N_j})} {(A_{j+1}q^{2N_j}, A_{j+2} q^{N_j + x_{j+2}};q)_x\o
(q, q^{N_j +1 - x_{j+2}}/a_{j+2};q)_x}\cr
&\times {(bA_{j+1}q^{2N_j+1}/a_1, q^{N_j - x_{j+2}};q)_x\o
(a_1/b, A_{j+1}q^{N_j + x_{j+2} +1};q)_x} \Big({bA_{j+2}\o a_1} 
q^{2N_j}\Big)^{-x}\cr
&\times r_{n_{j+1}} (x; bA_{j+1}q^{2N_j}/a_1, a_{j+2}q^{-1}, 
A_{j+1}q^{N_j + x_{j+2}},x_{j+2} - N_j;q)\cr
&\times r_{m_{j+1}} (x;bA_{j+1}q^{2N_j}/a_1, 
a_{j+2}q^{-1},A_{j+1}q^{N_j + x_{j+2}}, x_{j+2}-N_j;q)\cr
&= {(qA_{j+1}, qbA_{j+1}/a_1;q)_{2N_j} (a_{j+2};q)_{x_{j+2} - N_j} 
(A_{j+2};q)_{x_{j+2} + N_j}\o
(q;q)_{x_{j+2} - N_j} (qA_{j+1};q)_{x_{j+2} + N_j}}\cr
&\times\Big({bA_{j+1}\o a_1} q^{2N_j+1}\Big)^{-N_j} (q, bA_{j+1} 
q^{2N_j + 1}/a_1, a_{j+2}, A_{j+2}q^{N_j + x_{j+2}}, q^{N_j - 
x_{j+2}};q)_{n_{j+1}}\cr
&\times {(1-bA_{j+2}q^{2N_j}/a_1) (bA_{j+2}q^{N_j + x_{j+2} + 
1}/a_1;q)_{n_{j+1}}\o
(1-bA_{j+2}q^{N_j + x_{j+2}}/a_1) (bA_{j+2} q^{2N_j}/a_1;q)_{n_{j+1}}}\cr
&\times {(q^{-N_j - x_{j+2}}/A_{j+1}, bA_{j+2}q^{2N_j + 
1}/a_1;q)_{x_{j+2} - N_j}\o
(bq^{N_j - x_{j+2} + 1}/a_1, a_{j+2};q)_{x_{j+2} - N_j}} 
\delta_{n_{j+1}, m_{j+1}}\cr
}
$$
by (2.3). Substituting this into (2.12) and simplifying, we find that 
the right-hand side of (2.12) is (2.11) with $j$ replaced by $j+1$. 
This proves (2.11).
\medskip
Now set $j = s-1$ in (2.11) and substitute into the left-hand side of 
(2.9) to find that if we set $c_{s,s}= a_s$ for $s \ge 2$, then
$$\eqalignno{
&\sum_{\x} R_{\n} (\x;q) R_{\m}(\x;q) \rho (\x;q)\cr
&= \sum^N_{x_s = 0}\ \sum^{x_s}_{x_{s-1}=0} \ldots 
\sum^{x_3}_{x_2=0}\ \sum^{x_2}_{x_1 = 0}R_{\n}(\x;q) R_{\m} (\x;q) 
\rho (\x;q)\cr
&= \left[\prod^{s-1}_{k=1} {(1-bA_{k+1}/a_1) (q, a_{k+1};q)_{n_k} 
(qbA_k/a_1;q)_{N_k + N_{k-1}}\o
(1-bA_{k+1}q^{2N_k}/a_1) (bA_{k+1}/a_1;q)_{N_k + N_{k-1}}} 
\delta_{n_k, m_k}\right]\cr
&\times {(q,A_s, qbA_s/a_1;q)_{2N_{s-1}} (A_{s+1}q^N, q^{-N};q)_{N_{s-1}}\o
(A_s q^{N+1}, q^{1-N}/a_{s+1};q)_{N_{s-1}}} (A_{s+1}q^{N_{s-1}})^{-N_{s-1}}\cr
&\times \sum^{N-N_{s-1}}_{x=0} {(1-A_sq^{2N_{s-1} + 2x})\o
(1-A_s q^{2N_{s-1}})} {(A_s q^{2N_{s-1}}, {bA_s\o a_1} 
q^{2N_{s-1}+1}, A_{s+1} q^{N+N_{s-1}}, q^{N_{s-1}-N};q)_x\o
(q, a_1/b, q^{1-N + N_{s-1}}/a_{s+1}, A_s q^{N+N_{s-1}+1};q)_x}\cr
&\times \Big({bA_{s+1}\o a_1} q^{2N_{s-1}}\Big)^{-x} r_{n_s} (x; bA_s 
q^{2N_{s-1}}/a_1, a_{s+1}q^{-1}, A_s q^{N+N_{s-1}}, N-N_{s-1};q)\cr
&\times r_{m_s} (x; bA_s q^{2N_{s-1}}/a_1, a_{s+1}q^{-1}, 
A_sq^{N+N_{s-1}}, N-N_{s-1};q)\cr
&= \lm_{\n} (q) \delta_{\n, \m}&(2.13)\cr
}
$$
with
$$\eqalignno{
\lm_{\n} (q) &= \lm_{\n} (a_1, \ldots, a_{s+1}, b, N;q)\cr
&= {(qA_s, qbA_{s+1}/a_1;q)_N\o
(a_{s+1}, a_1/b;q)_N} \Big({a_1\o qbA_s}\Big)^N\cr
&\times (bA_{s+1}q^{N+1}/a_1, A_{s+1}q^N,bq^{1-N}/a_1, q^{-N};q)_{N_s}\cr
&\times\prod^s_{k=1} {(q, a_{k+1};q)_{n_k} (qbA_k/a_1;q)_{N_k + 
N_{k-1}} (1-bA_{k+1}/a_1)\o
(bA_{k+1}/a_1;q)_{N_k + N_{k-1}} (1-bA_{k+1}q^{2N_k}/a_1)},&(2.14)\cr
}
$$
where, in the weight function $\rho (\x;q)$ in (2.8) we set
$$c_{1,s} = bq \quad {\rm for}\ s \ge 1,\quad c_{k,s} = a_k\quad {\rm 
for}\ 2 \le k \le s, \eqno (2.15)$$
which completes the proof of (2.9) with
$$\eqalignno{
\rho (\x;q) &= {(q,qA_s;q)_N\o
(a_{s+1}, A_{s+1};q)_N} {(a_1, bq;q)_{x_1}\o
(q, a_1/b; q)_{x_1}} \Big({a_1\o bq}\Big)^{x_1}\cr
&\times \prod^s_{k=1} {(a_{k+1};q)_{x_{k+1} - x_k} 
(A_{k+1};q)_{x_{k+1} + x_k} (1-A_kq^{2x_k})\o
(q;q)_{x_{k+1} - x_k} (qA_k;q)_{x_{k+1} + x_k} (1-A_k)} a_k^{-x_k}.&(2.16)\cr
}
$$
\indent
Analogous to Tratnik's [24] permutations of the parameters and 
variables, we consider the following permutations
$$\eqalignno{
&a_1 \longleftrightarrow (A_s q^{2N})^{-1},\ \ \ a_{k+1} 
\longleftrightarrow a_{s-k+1},\ \  k = 1,
\ldots, s-1,\cr &a_{s+1} \longleftrightarrow bq,\ \ x_k 
\longleftrightarrow N - x_{s-k+1},
\ \  k = 1, \ldots,
s.&(2.17)\cr }
$$
Clearly the summation region in (2.9) remains unchanged under these 
permutations, and the weight function $\rho (\x;q)$ in (2.16) changes 
to a multiple of itself, namely,
$${(a_{s+1}, a_1/b, A_{s+1};q)_N\o
(a_1, bq, qA_s;q)_N} \Big({qbA_s\o a_1}\Big)^N \rho (\x; a_1, \ldots, 
a_{s+1}, b, N;q).\eqno (2.18)$$
On the other hand the polynomial $R_{\n}(\x;q)$ transforms to
$$R_{\n} ((N-x_s, N-x_{s-1}, \ldots, N-x_1); (A_sq^{2N})^{-1}, a_s, 
a_{s-1}, \ldots, a_2, bq, a_{s+1}q^{-1},N;q),\eqno (2.19)$$
which, since (2.6) is not invariant under (2.17), gives a second family 
of multivariable orthogonal
$q$-Racah polynomials. By replacing $n_k$ by $n_{s-k+1}$ for $k = 1, 
\ldots, s$, setting $N^*_k =
\prod\limits^s_{j=k} n_j$ and using (2.9) we obtain that the polynomials
$$\eqalignno{
&\tilde R_{\n} (\x;q) = \tilde R_{\n}(\x; a_1, \ldots, a_{s+1}, b, N;q)\cr
&= r_{n_1} (N-N^*_2 - x_1; A_{s+1}q^{2N^*_2-1}/a_1, b, 
q^{N^*_2-N}/a_1, N-N^*_2;q)\cr
&\times\prod^s_{k=2} r_{n_k} (N-N^*_{k+1} - x_{kj}; 
A_{s+1}q^{2N^*_{k+1}-1}/A_k, a_k q^{-1}, q^{N^*_{k+1}-N-x_{k-1}}/A_k, 
N - N^*_{k+1} - x_{k-1} ;q)\cr
&&(2.20)\cr
}
$$
are a $q$-extension of Tratnik's [24, (2.12)] second multivariable 
Racah polynomials and that they satisfy the orthogonality relation
$$\sum_{\x} \tilde R_{\n} (\x;q)\tilde R_{\m} (\x;q)\rho (\x;q)= 
\tilde\lm_{\n}(q) \delta_{\n,\m}\eqno (2.21)$$
for $N_s$, $M_s \le N$, where
$$\eqalignno{
&\tilde\lm_{\n}(q) = \tilde\lm_{\n} (a_1, \ldots, a_{s+1}, b, N;q)\cr
&= {(qA_s, qbA_{s+1}/a_1;q)_N\o
(a_1/b, a_{s+1};q)_N} \Big({a_1\o qbA_s}\Big)^N\cr
&\times (bA_{s+1}q^{N+1}/a_1, bq^{1-N}/a_1, A_{s+1}q^N, q^{-N};q)_{N^*_1}\cr
&\times {(q, qb;q)_{n_1} (A_{s+1}/a_1;q)_{N^*_1 + N^*_2} (1-b A_{s+1}/a_1)\o
(bA_{s+1}/a_1; q)_{N^*_1 + N^*_2} (1-bA_{s+1}q^{2N^*_1}/a_1)}\cr
&\times\prod^{s-1}_{k=1} {(q, a_{k+1};q)_{n_{k+1}} 
(A_{s+1}/A_{k+1};q)_{N^*_{k+1} + N^*_{k+2}}(1-A_{s+1}/qA_k)\o
(A_{s+1}/qA_k;q)_{N^*_{k+1} + 
N_{k+2}^*}(1-A_{s+1}q^{2N^*_{k+1}}/qA_k)}.&(2.22)\cr
}
$$
{\bf 3. Some limit cases of (2.9) and (2.21)}
\medskip
\noindent
The multivariable $q$-Racah polynomials in (2.6) and (2.20) contain 
as limit cases Tratnik's [24] systems of multivariable Racah, Hahn, 
dual Hahn, Krawtchouk, Meixner, and Charlier polynomials. Here we 
will derive limit cases of the orthogonality relations (2.9) and 
(2.21) containing multivariable $q$-Hahn, dual $q$-Hahn, 
$q$-Krawtchouk, $q$-Meixner, and $q$-Charlier polynomials. We start 
with the $q$-analogue of the dual Hahn polynomials defined by
$$\eqalignno{
&d_n (x;b, c, N;q) = \lim_{a\to 0} r_n(x; a, b, c, N;q)\cr
&= (bcq,q^{-N};q)_n \Big(q^N/c\Big)^{n/2} \ _3\phi_2\left[\matrix{
q^{-n}, q^{-x}, cq^{x-N}\cr
bcq,q^{-N}\cr} ;q,q\right]&(3.1)\cr
}
$$
for $n = 0, 1, \ldots, N$, and let
$$\eqalignno{
D_{\n} (\x;q) &= D_{\n} (\x; a_1, \ldots, a_{s+1}, N;q)\cr
&= \lim_{b\to 0} R_{\n} (\x; a_1, \ldots, a_{s+1}, b, N;q)\cr
&= \prod^s_{k=1} d_{n_k} (x_k - N_{k-1};a_{k+1} q^{-1}, A_k 
q^{x_{k+1} + N_{k-1}}, x_{k+1} - N_{k-1};q)&(3.2)\cr
}
$$
where $N_s \le N$ and $x_{s+1} = N$. Then, by taking the $b \to 0$ 
limit of (2.9) we get the orthogonality relation
$$\sum_{\x} D_{\n} (\x;q) D_{\m} (\x;q) \rho_D (\x;q) = \lm_D (\n;q) 
\delta_{\n, \m}\eqno (3.3)$$
for $N_s$, $M_s \le N$, were
$$\eqalignno{
\rho_D (\x;q) &= \rho_D (\x;a_1, \ldots, a_{s+1}, N;q)\cr
&= {(q, qA_s;q)_N\o
(a_{s+1}, A_{s+1};q)_N} {(a_1;q)_{x_1}\o
(q;q)_{x_1}} (-q)^{-x_1} q^{-{\ss x_1\mathstrut \choose\ss 2}}\cr
&\times \prod^s_{k=1} {(a_{k+1};q)_{x_{k+1} - x_k} 
(A_{k+1};q)_{x_{k+1} + x_k}(1-A_k q^{2x_k})\o
(q;q)_{x_{k+1} - x_k} (qA_k;q)_{x_{k+1} + x_k} (1-A_k)}a^{-x_k}_k,&(3.4)\cr
}
$$
$$\eqalignno{
\lm_D(\n;q) &= \lm_D(\n;a_1, \ldots, a_{s+1}, N;q)\cr
&= {(qA_s;q)_N (A_{s+1}q^N, q^{-N};q)_{N_s}\o
(a_{s+1};q)_N} (-qA_s)^{-N} q^{-{\ss N\mathstrut \choose\ss 2}}\cr
&\quad \times \prod^s_{k=1} (q, a_{k+1};q)_{n_k},&(3.5)\cr
}
$$
and the summation in (3.3) is over the same region as in (2.13).
\medskip
On the other hand, as $b\to 0$ the orthogonality relation (2.21) 
approaches the limit
$$\sum_{\x} T_{\n} (\x;q) T_{\m} (\x;q) \rho_T (\x;q) = \lm_T (\n;q) 
\delta_{\n,\m}\eqno (3.6)$$
for $N_s$, $M_s\le N$, where
$$\eqalignno{
T_{\n} (\x;q) &= T_{\n} (\x;a_1, \ldots, a_{s+1}, N;q)\cr
&= d_{n_1} (N-N^*_2 - x_1; A_{s+1}q^{N+N^*_2-1}, q^{N^*_2-N}/a_1, N-N^*_2;q)\cr
&\times \prod^s_{k=2} r_{n_k} (N-N^*_{k+1} - x_k; {A_{s+1}\o A_k} 
q^{2N^*_{k+1}-1}, {a_k\o q}, {q^{N^*_{k+1} - N - x_{k-1}}\o
A_k}, N-N^*_{k+1} - x_{k-1};q),\cr
&&(3.7)\cr
}
$$
$\rho_T(\x;q) = \rho_D (\x;q)$ is the same weight function as in (3.4), and
$$\eqalignno{
\lm_T (\n;q) &= \lm_T (\n; a_1, \ldots, a_{s+1}, N;q)\cr
&= {(qA_s;q)_N\o (a_{s+1};q)_N} (-qA_s)^{-N} q^{-{\ss N\mathstrut 
\choose\ss 2}}\cr
&\times (A_{s+1}q^N, q^{-N};q)_{N^*_1} (q;q)_{n_1} 
(A_{s+1}/a_1;q)_{N^*_1 + N^*_2}\cr
&\times \prod^{s-1}_{k=1} {(q, a_{k+1};q)_{n_{k+1}} 
(A_{s+1}/A_{k+1};q)_{N^*_{k+1} + N^*_{k+2}} (1-A_{s+1}/qA_k)\o
(A_{s+1}/qA_k;q)_{N^*_{k+1} + N_{k+2}^*} (1-A_{s+1} 
q^{2N^*_{k+1}}/qA_k)}.&(3.8)\cr
}
$$
\indent
If we multiply (2.9) by
$$\prod^s_{k=1} \Big({bA_k\o a_1} q^{N_{k-1}+1}\Big)^{-n_k} 
\Big({bA_k\o a_1}q^{M_{k-1} + 1}\Big)^{-m_k}$$
and take the limit $b \to \infty$, we obtain the limiting relation
$$\sum_{\x} D^*_{\n} (\x;q) D^*_{\m}(\x;q) \rho_{D^*}(\x,q) = 
\lm_{D^*}(\n;q) \delta_{\n,\m}\eqno (3.9)$$
for $N_s$, $M_s \le N$, where
$$\eqalignno{
D^*_{\n} (\x;q) &= D^*_{\n} (\x; a_1, \ldots, a_{s+1}, N;q)\cr
&= \prod^s_{k=1} d^*_{n_k} (x_k - N_{k-1}; a_{k+1} q^{-1}, A_k 
q^{x_{k+1} + N_{k-1}}, x_{k+1} - N_{k-1};q)&(3.10)\cr
}
$$
with
$$\eqalignno{
&d^*_n (x; b, c, N;q) = \lim_{a\to\infty} (aq)^{-n} r_n (x; a, b, c, N;q)\cr
&= (bcq, q^{-N};q)_n (-1)^n q^{{\ss n\mathstrut \choose\ss 2}} 
(q^N/c)^{n/2} \, _3\phi_2\left[\matrix{
q^{-n}, q^{-x}, cq^{x-N}\cr
bcq, q^{-N}\cr} ;q,bq^{n+1}\right]&(3.11)\cr
}
$$
for $n = 0, 1, \ldots, N$,
$$\eqalignno{
&\rho_{D^*} (\x;q) = \rho_{D^*}(\x; a_1, \ldots, a_{s+1}, N;q)\cr
&= {(q, qA_s;q)_N\o
(a_{s+1}, A_{s+1};q)_N} {(a_1;q)_{x_1}\o
(q;q)_{x_1}} (-a_1)^{x_1} q^{-{\ss x_1\mathstrut \choose\ss 2}}\cr
&\times \prod^s_{k=1} {(a_{k+1};q)_{x_{k+1} - x_k} 
(A_{k+1};q)_{x_{k+1} + x_k} (1-A_kq^{2x_k})\o
(q;q)_{x_{k+1} - x_k} (qA_k;q)_{x_{k+1} + x_k} (1-A_k)} a^{-x_k}_k, &(3.12)\cr
}
$$

and
$$\eqalignno{
&\lm_{D^*} (\n;q) = \lm_{D^*} (\n; a_1, \ldots, a_{s+1}, N;q)\cr
&= {(qA_s;q)_N\o (a_{s+1};q)_N} (-a_{s+1})^N q^{{\ss N\mathstrut 
\choose\ss 2}}A^{N_s}_{s+1} q^{N_s (N_s+1)}\cr
&\times \prod^s_{k=1} (q, a_{k+1};q)_{n_k} (A_kq^{N_{k-1} + 
1})^{-2n_k}q^{N_{k-1}-N_k} a_{k+1}^{-N_k - N_{k-1}}.&(3.13)\cr
}
$$
Note that $\rho_{D^*}(\x;q) = (a_1/q)^{x_1})$ $q^{{\ss x_1\mathstrut 
\choose\ss 2}}$ $\rho_D (\x;q)$ and that $D^*_\n(\x;q)$
is not a constant (independent of $\x$) multiple of $D_{\n}(\x;q)$. 
Similarly, by inverting the base
$q$ in the orthogonality relation (3.9) it follows that (3.9) is 
equivalent to (3.3).
\medskip
To obtain $q$-analogues of the multivariable Hahn polynomials in 
Karlin and McGregor [14] and in Tratnik [24], multiply (2.9) by 
$a_1^{N_s}$, take the limit $a_1\to 0$, replace $b$ by $a_1$ and 
$a_k$ by $qa_k$ for $k = 2, 3, \ldots, s+1$, and make the change of 
variables $y_1 = x_1$, $y_k = x_k - x_{k-1}$ for $k = 2, 3, \ldots, 
s$. This yields the multivariable $q$-Hahn polynomial orthogonality 
relation
$$\sum_{\bf y} H_{\n} ({\bf y};q) H_{\m} ({\bf y};q) \rho_H ({\bf 
y};q) = \lm_H (\n;q) \delta_{\n,\m}\eqno (3.14)$$
for $N_s$, $M_s \le N$, where
$$\eqalignno{
H_{\n}(\x;q) &= H_{\n}(\x;a_1, \ldots, a_{s+1}, N;q)\cr
&= \prod^s_{k=1} h_{n_k} (Y_k-N_{k-1}; A_kq^{2N_k + k-1}, a_{k+1}, 
Y_{k+1} - N_{k-1};q)&(3.15)\cr
}
$$
with the single-variable $q$-Hahn polynomials defined by
$$\eqalignno{
&h_n (x; a, b, N;q) = \lim_{c\to 0} (cq^{-N})^{n/2} r_n (x; a, b, c, N;q)\cr
&= (aq, q^{-N};q)_n \ _3\phi_2\left[\matrix{
q^{-n}, abq^{n+1}, q^{-x}\cr
aq, q^{-N}\cr} ;q,q\right]&(3.16)\cr
}
$$
for $n = 0, 1, \ldots, N$,
$$\eqalignno{
\rho_H ({\bf y};q) &= \rho_H ({\bf y}; a_1, \ldots, a_{s+1}, N;q)\cr
&= {(q^{-N};q)_{Y_s}\o (q^{-N}/a_{s+1};q)_{Y_s}} a^{-Y_s}_{s+1} 
\prod^s_{k=1} {(qa_k;q)_{y_k}\o
(q;q)_{y_k}} (qa_k)^{-Y_k},&(3.17)\cr
}
$$
$$\eqalignno{
\lm_H (\n;q)&= \lm_H(\n; a_1, \ldots, a_{s+1}, N;q)\cr
&= {(A_{s+1}q^{s+1};q)_{N+N_s}(q^{-N};q)_{N_s}\o
(qa_{s+1};q)_N} (A_s q^{N_s + s)})^{-N}(-1)^{N_s} q^{{\ss 
N_s\mathstrut \choose\ss 2}}\cr
&\times\prod^s_{k=1} {(q,qa_{k+1};q)_{n_k} (A_kq^k;q)_{N_k + N_{k-1}} 
(1-A_{k+1}q^k)\o
(A_{k+1}q^k;q)_{N_k + N_{k-1}} (1-A_{k+1}q^{k+2N_k})} 
\Big(A_kq^{k+2N_{k-1}}\Big)^{n_k},&(3.18)\cr
}
$$
$Y_k = \sum\limits^k_{j=1} y_j$ for $1 \le k \le s$ so that $Y_1 = y_1 = x_1$
and $x_k = Y_k$ for $1 \le k \le s$, $Y_{s+1} = x_{s+1} = N$, and the summation is 
over all ${\bf y}$
with $y_k = 0, 1, \ldots$, and $Y_s \le N$. Also see Dunkl [4] and 
Rosengren [18].
\medskip
Clearly the summation region in (3.14) and the weight function in 
(3.17) are invariant under any permutation of the labels $(1,2, 
\ldots, s)$. If we set $y_{s+1} = N-Y_s$, then they are also 
invariant (apart from the renormalization of the weight function) 
under any permutation of the labels $(1,2, \ldots, s+1)$; i.e., under 
any simultaneous permutation of $(a_1, a_2, \ldots, a_{s+1})$ and 
$(y_1, y_2, \ldots, y_{s+1})$. Since the polynomials in (3.15) are 
generally not invariant under these permutations, this generates 
distinct systems satisfying the orthogonality relation (3.14) via 
these permutations, which are $q$-extensions of those mentioned at 
the top of page 2341 in Tratnik [24].
\medskip
Next, consider the $q$-Krawtchouk polynomials defined by
$$\eqalignno{
&k_n (x; b, N;q) = \lim_{a\to\infty} (aq)^{-n} h_n (x; a, b, N;q)\cr
&= (q^{-N};q)_n (-1)^n q^{{\ss n\mathstrut \choose\ss 2}}\ 
_2\phi_1\Big(q^{-n}, q^{-x};q^{-N};q, bq^{n+1}\Big)&(3.19)
\cr}
$$
for $n = 0, 1, \ldots, N$. In view of (3.19), if we multiply (3.14) 
by $a^{-2N_s}_1$, let $a_1\to \infty$, and then replace $a_k$ by 
$a_{k-1}$ for $2\le k \le s+1$, we obtain the multivariable 
$q$-Krawtchouk orthogonality relation
$$\sum_{\bf y} K_{\n} ({\bf y};q) K_{\m}({\bf y};q) \rho_K ({\bf 
y};q) = \lm_K (\n;q)\delta_{\n,\m}\eqno (3.20)$$
for $N_s$, $M_s \le N$, where
$$\eqalignno{
K_{\n} ({\bf y};q) &= K_{\n}({\bf y}; a_1, \ldots, a_s, N;q)\cr
&= \prod^s_{j=1} k_{n_j} (Y_j - N_{j-1}; a_j, Y_{j+1} - N_{j-1};q),&(3.21)\cr
}
$$
$$\eqalignno{
\rho_K ({\bf y};q) &= \rho_K ({\bf y}; a_1, \ldots, a_s, N;q)\cr
&= {(q^{-N};q)_{Y_s}\o
(q;q)_{y_1} (q^{-N}/a_s;q)_{Y_s}} (-1)^{y_1} q^{-{\ss y_1\mathstrut 
\choose\ss 2}}a_s^{-Y_s}\cr
&\quad \times \prod^s_{j=2} {(qa_{j-1};q)_{y_j}\o
(q;q)_{y_j}} (qa_{j-1})^{-Y_j},&(3.22)\cr
}
$$
$$\eqalignno{
\lm_K (\n;q) &= \lm_K(\n; a_1, \ldots, a_s, N;q)\cr
&= {(q^{-N};q)_{N_s}\o
(qa_s;q)_{N_s}} (A_{s-1} q^{s+N_s})^{-N} (A_sq^{s+1})^{N+N_s} (-1)^N 
q^{{\ss N_s\mathstrut \choose\ss 2} + {\ss N+N_s\mathstrut \choose 
\ss 2}}\cr
&\quad\times\prod^p_{j=1} (q, qa_j;q)_{n_j} 
(A_{j-1}q^{j+2N_{j-1}})^{-n_j} a_j^{-N_j - N_{j-1}} 
q^{-2N_j}&(3.23)\cr
}
$$
and the summation is over the same region as in (3.14). This 
orthogonality relation can also be derived by starting with the $a_1 
\to 0$ limit case of (3.14).
\medskip
Now consider the $q$-Meixner polynomials defined by
$$M_n (q^{-x}; a, c;q)=\, _2\phi_1 \Big(q^{-n}, q^{-x};aq; q, 
-q^{n+1}/c\Big)\eqno (3.24)$$
for $n = 0, 1, \ldots$, which satisfy the orthogonality relation
$$\eqalignno{
&\sum^\infty_{x=0} M_n (q^{-x}; a, c;q) M_m (q^{-x};a, c;q) 
{(aq;q)_x\o (q, -acq;q)_x} c^x q^{{\ss x\mathstrut \choose\ss 2}}\cr
&= {(-c;q)_\infty (q, -q/c;q)_n\o
(-acq;q)_\infty (aq;q)_n} q^{-n} \delta_{n,m}.&(3.25)\cr
}
$$
See [9, Ex. 7.12] and [15, \S3.13].
\medskip
To find multivariable orthogonal $q$-Meixner polynomials we, formally,  replace
the upper limit of summation $N$ for the sum over $x_s$ in (2.9) by $\infty$
and  replace $q^{-N}$ by an arbitrary number $\beta$, say. 
Then, from (2.16), the factor
$$\eqalignno{
&{(q, qA_s;q)_N\o
(a_{s+1}, A_{s+1};q)_N} {(a_{s+1};q)_{N-x_s} (A_{s+1};q)_{N+x_s}\o
(q;q)_{N-x_s} (qA_s;q)_{N+x_s}}\cr
&= {(q^{-N}, A_{s+1}q^N;q)_{x_s}\o
(q^{1-N}/a_{s+1}, A_sq^{N+1};q)_{x_s}} (q/a_{s+1})^{x_s},\cr
}
$$
gets replaced by
$${(\beta, A_{s+1}/\beta;q)_{x_s}\o
(\beta q/a_{s+1}, qA_s/\beta;q)_{x_s}} (q/a_{s+1})^{x_s},$$
and, from (2.14), the factor
$$\eqalign{
&{(qA_s, qbA_{s+1}/a_1;q)_N\o
(a_{s+1}, a_1/b;q)_N} \Big({a_1\o qbA_s}\Big)^N\cr
&= {(qA_s, a_1q^{-N}/b A_{s+1};q)_N\o
(q^{1-N}/a_{s+1}, a_1/b;q)_N}\cr
}
$$
gets replaced by
$${(qA_s, a_1\beta/b A_{s+1}, q/a_{s+1}, a_1/b\beta ;q)_\infty\o
(\beta q/a_{s+1}, a_1/b, qA_s/\beta, a_1/bA_{s+1};q)_\infty}$$
via the identity $(a;q)_N = (a;q)_\infty/(aq^N;q)_\infty$ in [9, 
(I.5)]. Hence, by multiplying both sides of (2.9) by 
$(bq/a_1^{1/2})^{-N_s - M_s}$, taking the limits $b \to \infty$ and 
$a_1 \to 0$, and replacing $a_k$ by $a_{k-1}$ for $2 \le k \le s+1$, 
we are led  to the limit case
$$\eqalignno{
&\sum^\infty_{x_s = 0} \sum^{x_s}_{x_{s-1} = 0} \ldots 
\sum^{x_2}_{x_1 = 0} {\cal M}_{\n} (\x; q) {\cal M}_{\m} (\x;q) 
\rho_{\cal M} (\x;q)\cr
&= \lm_{\cal M} (\n;q) \delta_{\n,\m}\,,&(3.26)\cr
}
$$
where
$$\eqalignno{
{\cal M}_{\n} (\x;q) &= {\cal M}_{\n} (\x;a_1, \ldots, a_s, \beta;q)\cr
&= \left[\prod^{s-1}_{k=1} (q^{N_{k-1} - x_{k+1}};q)_{n_k} (-1)^{n_k} 
q^{{\ss n_k\mathstrut \choose\ss 2} + n_k 
N_{k-1}}A_{k-1}^{n_k/2}\right.\cr
&\quad \left.\times M_{n_k} (q^{N_{k-1} - x_k};q^{N_{k-1}-x_{k+1} - 
1}, -q/a_k ;q)\right]\cr
&\times (\beta q^{N_{s-1}};q)_{n_s} (-1)^{n_s} q^{{\ss n_s\mathstrut 
\choose\ss 2} + n_s N_{s-1}} A_{s-1}^{n_s/2}\cr
&\times M_{n_s} (q^{N_{s-1} - x_s};\beta q^{N_{s-1}-1}, -q/a_s;q),&(3.27)\cr
}
$$
$$\eqalignno{
\rho_{\cal M} (\x;q) &= \rho_{\cal M}(\x;a_1, \ldots, a_s, \beta;q)\cr
&= {(\beta;q)_{x_s}\o
(q\beta/a_s;q)_{x_s}} \Big({q\o a_{s-1}a_s}\Big)^{x_s} {(-1)^{x_1}\o 
(q;q)_{x_1}}q^{{\ss x_1\mathstrut \choose\ss 2}}\cr
&\times \prod^{s-1}_{k=1} {(a_k;q)_{x_{k+1} - x_k}\o
(q;q)_{x_{k+1} - x_k}} a^{-x_k}_{k-1} &(3.28)\cr
}
$$
with $a_0 = 1$, and
$$\eqalignno{
\lm_{\cal M} (\n;q) &= \lm_{\cal M} (\n;a_1, \ldots, a_s, \beta;q)\cr
&= {(q/a_s;q)_\infty\o
(q\beta/a_s;q)_\infty} (\beta;q)_{N_s} q^{N_s(N_s-1)}A_s^{N_s}\cr
&\times \prod^s_{k=1} (q, a_k;q)_{n_k} a_k^{-N_k - 
N_{k-1}}q^{N_{k-1}-N_k}.&(3.29)\cr
}
$$
\noindent
This orthogonality relation can be verified by starting with the $q$-Meixner
orthogonality relation (3.25) and proceeding as in the derivation of (2.9).

The $\beta \to 0$ limit of (3.26) gives the following multivariable 
extension of the $q$-Charlier polynomial orthogonality
$$\eqalignno{
&\sum^\infty_{x_s = 0} \sum^{x_s}_{x_{s-1}=0} \ldots \sum^{x_2}_{x_1 
= 0} {\cal C}_{\n}(\x;q) {\cal C}_{\m}(\x;q)\rho_{\cal C}(\x;q)\cr
&= \lm_{\cal C}(\n;q)\delta_{\n,\m}\,,&(3.30)\cr
}
$$
where
$$\eqalignno{
{\cal C}_{\n} (\x;q) &= {\cal C}_{\n} (\x; a_1, \ldots, a_s;q)\cr
&= \left[\prod^{s-1}_{k=1} (q^{N_{k-1} - x_{k+1}};q)_{n_k} (-1)^{n_k} 
q^{{\ss n_k\mathstrut \choose\ss 2} + n_k N_{k-1}} 
A_{k-1}^{n_k/2}\right.\cr
&\quad\times\left. M_{n_k} (q^{N_{k-1} - x_k};q^{N_{k-1} - x_{k+1} - 
1}, -q/a_k;q)\right]\cr
&\times (-1)^{n_s} q^{{\ss n_s\mathstrut \choose\ss 2} + n_s 
N_{s-1}}A^{n_s/2}_{s-1}\cr
&\times c_{n_s} (q^{N_{s-1} - x_s}; -q/a_s;q),&(3.31)\cr
}
$$
$$\eqalignno{
\rho_{\cal C} (\x;q) &= \rho_{\cal C} (\x;a_1, \ldots, a_s;q)\cr
&=\Big({q\o a_{s-1} a_s}\Big)^{x_s} {(-1)^{x_1}\o (q;q)_{x_1}} 
q^{{\ss x_1\mathstrut \choose\ss 2}}\prod^{s-1}_{k=1} 
{(a_k;q)_{x_{k+1} - x_k}\o
(q;q)_{x_{k+1} - x_k}} a^{-x_k}_{k-1},&(3.32)\cr
}
$$
$$\eqalignno{
\lm_{\cal C} (\n;q) &=\lm_{\cal C} (\n;a_1, \ldots, a_s;q)\cr
&= (q/a_s;q)_\infty q^{N_s(N_s-1)}A_s^{N_s} \prod^s_{k=1} (q, 
a_k;q)_{n_k} a^{-N_k-N_{k-1}}_k q^{N_{k-1}-N_k},&(3.33)\cr
}
$$
with $a_0 = 1$ and the $q$-Charlier polynomial defined by
$$\eqalignno{
c_n(q^{-x}; a;q)&= \lim_{\beta \to 0} M_n (q^{-x}; \beta, a;q)\cr
&=\, _2\phi_1 (q^{-n}, q^{-x};0;q, -q^{n+1}/a).&(3.34)\cr
}
$$
See [9, Ex. 7.13] and [15, \S3.23] for the orthogonality relation and 
other properties of the $q$-Charlier polynomials. Of course 
(3.3), (3.6), (3.9), (3.14), (3.20), and (3.30) 
 can also be derived directly from their 
single variable special cases by proceeding as in the derivation of 
(2.9) and (2.21).
\medskip
Other discrete multivariable extensions of the Racah and Hahn 
polynomials are considered in [3], [5], [13], [17], [18], [20], [22] 
and [26]. For some related results see, e.g., [6], [7], [8], [9], 
[12], [16] and [19]. In view of the nonnegativity results for the 
linearization coefficients and for kernels containing products of 
$q$-Racah polynomials in [7] and [8] and the resulting convolution 
structure and positive summability methods, it would be of interest 
to see if any of these nonnegativity results can be extended to some 
of the multivariable orthogonal polynomials considered in this paper.
\bigskip
\noindent
{\bf References}
\medskip
\item{1.} R. Askey and J.A. Wilson, ``A set of orthogonal polynomials
that generalize the Racah
coefficients or $6$-$j$ symbols,''  {\it SIAM J. Math. Anal.}  {\bf
10} (1979), 1008--1016.
\medskip
\item{2.} R. Askey and J.A. Wilson,  ``Some basic hypergeometric
orthogonal polynomials that generalize
Jacobi polynomials,''  {\it Memoirs Amer. Math. Soc.}  {\bf 319} (1985).
\medskip
\item{3.} J.F. van Diejen and J.V. Stokman,  ``Multivariable $q$-Racah
polynomials,'' {\it Duke Math.\ J.} {\bf 91}  (1998), 89--136.
\medskip
\item{4.}  C.F. Dunkl, ``Orthogonal polynomials in two variables of
$q$-Hahn and $q$-Jacobi type,'' {\it SIAM J. Alg. Disc. Meth.} {\bf
1} (1980), 137--151.
\medskip
\item{5.}  C.F. Dunkl and Yuan Xu, {\it Orthogonal Polynomials of
Several Variables,}
Cambridge Univ. Press, 2001.
\medskip
\item{6.} G. Gasper,  ``Positivity and special functions,'' {\it
Theory and Applications of Special Functions} (R.\ Askey, ed.), Academic
Press, New York, 1975, 375--433.
\medskip
\item{7.}  G. Gasper and M. Rahman,  ``Nonnegative kernels in product
formulas for $q$-Racah polynomials,'' {\it J.\ Math.\ Anal.\ Appl.}
{\bf 95} (1983), 304--318.
\medskip
\item{8.} G. Gasper and M. Rahman, ``Product formulas of Watson,
Bailey and Bateman types and positivity of the Poisson kernel for
$q$-Racah polynomials,'' {\it SIAM J.\ Math.\ Anal.} {\bf 15} (1984), 768--789.
\medskip
\item{9.} G. Gasper and M. Rahman,  {\it Basic Hypergeometric
Series}\  (2nd Edition),
Cambridge Univ. Press, 2004.
\medskip
\item{10.} G. Gasper and M. Rahman,  ``$q$-Analogues of some multivariable
biorthogonal polynomials,''  {\it Theory and Applications of
Special  Functions.\ A volume dedicated to Mizan Rahman} (M.\ E.\ H.\
Ismail and E.\ Koelink,  eds.),
   {\it Dev.\ Math.} (2004), 217--240.
\medskip
\item{11.} G. Gasper and M. Rahman, ``Some systems of multivariable orthogonal
Askey--Wilson polynomials,''  {\it Theory and Applications of
Special  Functions.\ A volume dedicated to Mizan Rahman} (M.\ E.\ H.\
Ismail and E.\ Koelink,  eds.),
   {\it Dev.\ Math.} (2004), 331--341.
\medskip
\item{12.} Ya. I. Granovskii, and A.S. Zhedanov,
    `` `Twisted' Clebsch-Gordan coefficients for $SU_q(2)$,''
{\it J. Phys. A} {\bf 25} (1992), L1029--L1032.
\medskip
\item{13.} R.A. Gustafson,  ``A Whipple's transformation for
hypergeometric series in $U(n)$ and multivariable hypergeometric
orthogonal polynomials,'' {\it SIAM J.\ Math.\ Anal.} {\bf 18}
(1987), 495--530.
\medskip
\item{14.} S. Karlin and J. McGregor,  ``Linear growth models with
many types and
multidimensional Hahn polynomials,''
{\it Theory and Applications of Special Functions} (R.\ Askey, ed.), Academic
Press, New York, 1975, 261--288.
\medskip
\item{15.} R. Koekoek and R.F. Swarttouw, {\it  The Askey--scheme of
hypergeometric orthogonal
polynomials and its  $q$-analogue,}  Report 98-17, Delft Univ. of Technology,
   http://aw.twi.tudelft.nl/\~{} koekoek/research.html, 1998.
\medskip
\item{16.} H.T. Koelink and J. Van der Jeugt, ``Convolutions for orthogonal
polynomials from Lie and quantum algebra representations,'' {\it SIAM
J.\ Math.\ Anal.} {\bf 29} (1998), 794-822.
\medskip
\item{17.} T.H. Koornwinder, ``Askey--Wilson polynomials for root
systems of type $BC$,''
{\it Contemp. Math.} {\bf 138} (1998), 189--204.
\medskip
\item{18.} H. Rosengren,  ``Multivariable $q$-Hahn polynomials as coupling
coefficients for quantum algebra representations,'' {\it Int.\ J.\
Math.\ Math.\ Sci.}
{\bf 28} (2001), 331--358.
\medskip
\item{19.} V.P. Spiridonov,   ``Theta hypergeometric integrals,''
{\it Algebra i Analiz} ({\it St.\ Petersburg Math.\ J.}) {\bf 15}
(2003), 161--215.
\medskip
\item{20.} J.V. Stokman,  ``On $BC$ type basic hypergeometric orthogonal
polynomials,'' {\it Trans.\ Amer.\ Math.\ Soc.} {\bf 352} (2000), 1527--1579.
\medskip
\item{21.} M.V. Tratnik, ``Multivariable Wilson polynomials,''
   {\it J. Math. Phys.}  {\bf 30} (1989), 2001--2011.
\medskip
\item{22.} M.V. Tratnik, ``Multivariable biorthogonal
continuous---discrete Wilson and Racah polynomials,''
   {\it J. Math. Phys.}  {\bf 31} (1990), 1559--1575.
\medskip
\item{23.} M.V. Tratnik, ``Some multivariable orthogonal polynomials
of the Askey tableau---continuous
families,'' {\it J. Math. Phys.}  {\bf 32} (1991), 2065--2073.
\medskip
\item{24.} M.V. Tratnik,  ``Some multivariable orthogonal polynomials
of the Askey
tableau---discrete families,''
   {\it J. Math. Phys.}  {\bf 32} (1991), 2337--2342.
\medskip
\item{25.} J.A. Wilson, ``Some hypergeometric orthogonal
polynomials,''  {\it SIAM J. Math. Anal.}
   {\bf 11} (1980), 690--701.
\medskip
\item{26.} Yuan Xu, ``On discrete orthogonal polynomials of several
variables,'' to appear.

\bye